\documentclass[12pt]{article}
\usepackage{amssymb}
\usepackage{amsmath,amsthm}
\usepackage[latin1]{inputenc}

\usepackage{hyperref}
\hypersetup{colorlinks=true, urlcolor=blue}

 \newtheorem{remark}{Remark}

 \newtheorem{theorem}[remark]{Theorem}
 \newtheorem{proposition}[remark]{Proposition}
 \newtheorem{corollary}[remark]{Corollary}

\title{Offensive $k$-alliances in graphs}

\author{Henning Fernau$^1$, Juan A. Rodr\'{\i}guez$^2$ and Jos\'{e} M. Sigarreta$^3$\\
\\
$^1${\em FB 4-Abteilung Informatik}\\ Universit\"{a}t Trier,\\ 54286 Trier, Germany. \\
{\small e-mail:\mbox{\tt fernau\@@uni-trier.de}}
\\
$^2${\em Departament d'Enginyeria Inform\`{a}tica i Matem\`{a}tiques
}\\
Universitat Rovira i Virgili, \\ Av. Pa\"{\i}sos Catalans 26, 43007
Tarragona, Spain.
\\
{\small e-mail:\mbox{\tt juanalberto.rodriguez\@@urv.cat}}
\\
$^3${\em Departamento de Matem\'{a}ticas}\\ Universidad Carlos III
de Madrid,\\ Avda. de la Universidad 30, 28911 Legan\'es (Madrid),
Spain.
\\
 {\small e-mail:\mbox{\tt josemaria.sigarreta\@@uc3m.es}}
}

\date{}

\begin{document}

\maketitle

\begin{abstract}
Let   $G=(V,E)$ be a simple graph.  For a nonempty set $X\subset V,$
and a vertex $v\in V,$ $\delta_{X}(v)$ denotes the number of
  neighbors  $v$ has in $X.$
A nonempty set $S\subset V$ is an \emph{offensive  $r$-alliance} in
$G$ if $\delta _S(v)\ge \delta_{\bar{S}}(v)+r,$ $\forall v\in
\partial (S),$ where $\partial (S)$ denotes the boundary of $S$. An
offensive $r$-alliance $S$ is called \emph{global} if it forms a
dominating set. The \emph{global offensive $r$-alliance number} of
$G$, denoted by $\gamma_{r}^{o}(G)$, is the minimum cardinality of a
global offensive $r$-alliance in $G$. We show that the problem of
finding optimal (global) offensive $r$-alliances is NP-complete and we obtain
several tight bounds on $\gamma_{r}^{o}(G)$.
\end{abstract}

{\it Keywords:} Computational complexity, offensive alliances,
alliances in graphs, domination.

{\it AMS Subject Classification Numbers:} 03D15;  05C69; 05A20

\section{Introduction}

The mathematical properties of alliances in graphs were first
studied by P. Kristiansen, S. M. Hedetniemi and S. T. Hedetniemi
\cite{alliancesOne}. They proposed diffe\-rent types of alliances:
namely, defensive alliances
\cite{note,GlobalalliancesOne,alliancesOne,albesiga8}, offensive
alliances \cite{chellali,chellali1,fava,albesiga6,albesiga7} and
dual alliances or po\-werful alliances \cite{cota}. A generalization
of these alliances called $r$-alliances was presented by K. H.
Shafique and R. D. Dutton \cite{kdaf,kdaf1}.

In this paper, we  study the mathematical properties of offensive
$r$-alliances. We begin by stating the  terminology used. Throughout
this article, $G=(V,E)$ denotes a simple graph of order $|V|=n$. We
denote two adjacent vertices $u$ and $v$ by $u\sim v$. For a
nonempty set $X\subseteq V$, and a vertex $v\in V$,
 $N_X(v)$ denotes the set of neighbors  $v$ has in $X$:
$N_X(v):=\{u\in X: u\sim v\},$ and the degree of $v$ in $ X$ will be
denoted by $\delta_{X}(v)=|N_{X}(v)|.$ We denote the degree of a
vertex $v\in V$  by $\delta(v)$ and the degree sequence of
 $G$ by $\delta_{1}\geq
\delta_{2}\geq \cdots \geq\delta_{n}$. The complement of the
vertex-set $S$ in $V$ is denoted by $\bar{S}$ and the boundary,
$\partial (S)$, of $S$ is defined by
$$\partial(S):=\displaystyle{\bigcup_{v\in S}} N_{\bar{S}}(v).$$

For $r\in \{2-\delta_1,\dots,\delta_1\}$, a nonempty set
{ $S\subset V$ is an \emph{offensive  $r$-alliance}
in $G$   if for every $ v\in
\partial (S)$,}
\begin{equation}\label{cond-A-Defensiva} \delta _S(v)\ge \delta_{\bar{S}}(v)+r.\end{equation}
 or, equivalently,
\begin{equation}\label{cond-A-Defensiva1} \delta (v)\ge 2\delta_{\bar{S}}(v)+r.\end{equation}
An offensive $1$-alliance is an \emph{offensive alliance} and an
offensive $2$-alliance is a \emph{strong offensive alliance} as
defined in \cite{fava,albesiga6,albesiga7}.

The \emph{offensive $r$-alliance number} of $G$, denoted by
$a_r^o(G)$, is defined as the minimum cardinality of an offensive
$r$-alliance in $G$. Notice that
\begin{equation}a^o_{r+1}(G)\ge a^o_r(G).\end{equation}
The offensive $1$-alliance number of $G$ is known as the
\emph{offensive alliance number} of $G$
 and the offensive  $2$-alliance number is known as  the \emph{strong offensive alliance number}
  \cite{fava,albesiga6,albesiga7}.

A set $S\subset V$ is a  \emph{dominating
set}\label{conjuntodominante} in $G=(V,E)$ if for every vertex $u\in
\bar{S}$,  $\delta_S(u)>0$ (every vertex in $\bar{S}$ is adjacent to
at least one vertex in S). The \emph{domination number} of $G$,
denoted by $\gamma(G)$, is the minimum cardinality of a dominating
set in $G$.

An offensive $r$-alliance $S$ is called \emph{global} if it forms a
dominating set, i.e., $\partial(S)=\bar{S}$. The \emph{global
offensive $r$-alliance number} of $G$, denoted by
$\gamma_{r}^{o}(G)$, is the minimum cardinality of a global
offensive $r$-alliance in $G$. Clearly,
\begin{equation}
\gamma_{r+1}^o(G)\ge \gamma_r^o(G)\ge \gamma(G)\quad {\rm and }
\quad \gamma_r^o(G)\ge a_r^o(G).\end{equation}

Notice that if every vertex of $G$ has even
 degree and $k$ is odd, $k=2l-1$ ,
 then
every offensive $(2l-1)$-alliance in $G$ is an offensive
$(2l)$-alliance. Hence, in such a case, $a_{2l-1}^o(G)=a_{2l}(G)$
and $\gamma^{o}_{2l-1}(G)=\gamma^{o}_{2l}(G).$ Analogously, if every
vertex of $G$ has odd
 degree and $k$ is even, $k=2l$,
 then every offensive $(2l)$-alliance in $G$ is an offensive
$(2l+1)$-alliance. Hence, in such a case,
$a_{2l}^o(G)=a_{2l+1}^o(G)$ and
$\gamma^{o}_{2l}(G)=\gamma^{o}_{2l+1}(G).$

\section{On the complexity of finding optimal offensive $r$-alliances}

For the class of complete graphs of order $n$, $G=K_{n}$, we have
the exact value of $a_r^o(G)$. That is,
\begin{align*}
n-1=&a^{o}_{n-1}(K_n)=a^{o}_{n-2}(K_n) \\
 \ge&a^{o}_{n-3}(K_n)=a^{o}_{n-4}(K_n)=n-2\\
&\cdots \\
\ge&a^{o}_{5-n}(K_n)=a^{o}_{4-n}(K_n)=2\\
\ge &a^{o}_{3-n}(K_n)=1.
\end{align*}
Hence, for every $r\in \{3-n, \dots, n-1\}$,
 $a_{r}^{o}(K_n)=\left\lceil\frac{n+r-1}{2}\right\rceil.$
 In this case, every offensive $r$-alliance is global and every vertex-set
  of cardinality $\left\lceil\frac{n+r-1}{2}\right\rceil$ is a (global)
 offensive $r$-alliance.

As we will se below, in general, the problem of finding optimal
(global) offensive $r$-alliances is NP-complete. That is, we are
interested in the computational complexity of the following
optimization problems.

\medskip

\textbf{Offensive $r$-Alliance problem ($r$-OA)}:
\begin{itemize}
\item[]Given: A graph $G=(V,E)$ and a positive integer $k\le |V|$.
\item[]Question: Is there an offensive $r$-alliance in $G$ of size $k$ or
less?
\end{itemize}

\textbf{Global offensive $r$-Alliance problem ($r$-GOA)}:
\begin{itemize}
\item[]Given: A graph $G=(V,E)$ and a positive integer $k\le |V|$.
\item[]Question: Is there a global offensive $r$-alliance in $G$ of size $k$ or
less?
\end{itemize}

\subsection{Offensive alliances}

Our reasoning will use and generalize the following observation:

\begin{proposition}{\rm \cite{fava}}\label{prop2}
On cubic graphs, every vertex cover is a strong offensive alliance
and vice versa.
\end{proposition}

With some gadgetry, this was used in~\cite{FerRai07} to show
NP-hardness of finding small offensive alliances. We will generalize
those results in the following.

\begin{theorem}
\label{offensive-NP} $\forall r$: $r$-OA is NP-complete.
\end{theorem}

\begin{proof}
It is clear that $r$-OA is in NP.

Consider first the case that $r\geq 3$. For any connected
$r$-regular graph $G=(V,E)$, it can be seen that $C\subseteq V$ is a
minimum vertex cover iff $C$ is a minimum $r$-offensive alliance.
Clearly, any vertex cover is an $r$-OA. Let $S$ be an $r$-OA. By
definition, $S\neq \emptyset$. Discuss $x\in S$. Any neighbor of $x$
must have $r$, i.e., all, neighbors in $S$, and we can continue the
argument with those vertices taking the role of $x$, till the whole
graph is exhausted (since it is connected by assumption). Hence, the
complement of $S$ forms an independent set, which means that $S$
itself is a vertex cover. Since it is well-known that the vertex
cover problem, restricted to $r$-regular graphs is NP-complete for any $r\geq 3$,
see \cite{Fei2003} for a recent account related to approximability results,
the claim follows for $r\geq 3$.

Now, we show: if $r$-OA is NP-hard, then so is $(r-1)$-OA. By
induction, the whole claim will follow.

Let $(G=(V,E),k)$ be an instance of $r$-OA, with $n=|V|$. We
construct an instance of $(r-1)$-OA as follows: $G'=(V',E')$ with
$V'=V\times\{1,2,3\}\cup \{c_1,\dots,c_{n-r+2}\}$.
In $E'$, we find the following edges (and only those):\\
--- $\{(u,1),(v,1)\}\in E'$ iff $\{(u,2),(v,2)\}\in E'$ iff $\{u,v\}\in E$;\\
--- $\{(u,1),(u,3)\}\in E'$ and $\{(u,2),(u,3)\}\in E'$ for any $u\in V$;\\
--- $\{(u,3),c_j\}\in E'$ for any $u\in V$ and any $1\leq j\leq n-r+2$;\\
--- $\{c_i,c_j\}\in E'$  for any $1\leq i< j\leq n-r+2$.\\
Let $k'=2k$. As  in~\cite{FerRai07}, one can show that $S$ is an
$r$-OA of size at most $k$ for $G$ iff $S\times\{1,2\}$ is a $(r-1)$-OA of
size at most $2k$ for $G'$, and that there is no other possibility to form
smaller $(r-1)$-OAs in $G'$ due to the attached clique.
\end{proof}

\subsection{Global offensive alliances}

Cami et al.~\cite{Cametal2006} showed NP-completeness for $r=1$. We
are going to modify their construction to show NP-completeness for
any fixed $r$. Since we are dealing with the degree of vertices both
in $G$ and within the new graph $G'$ as constructed below, we are
going to attach $G$ and $G'$ to $\delta$ to avoid confusion in our
notation.

\begin{theorem}
\label{globaloffensive-NP} $\forall r$: $r$-GOA is NP-complete.
\end{theorem}

\begin{proof}
Membership in NP is clear.

The construction in~\cite{Cametal2006}
 can
be 
modified to work for any case $r\leq 1$. Let $(G,k)$ be an
instance of Dominating Set with minimum degree $|r|+1$, with
$G=(V,E)$. To any $v\in V$, attach $\delta_G(v)+r-1\geq 0$ copies of
$K_2$ with one edge per $K_2$-copy, this way yielding a new graph $G'=(V',E')$ with $G$ as a
subgraph; call the new neighbors of vertices from $V$ $A$-vertices
and collect them into set $A$, and call $N(A)\setminus V$
$B$-vertices.

If $D\subseteq V$ is a dominating set in $G$, then $S=D\cup A$ is a
$r$-GOA. Clearly, $S$ is a dominating set in $G'$. Now, consider a
$B$-vertex $v$. Obviously, $N(v)\subseteq A$, and therefore
$|N_{G'}(v)\cap S|\geq |N_{G'}(v)\cap\bar S|+r$. Any vertex  $v\in
V\setminus D$ has a neighbor $d\in D$. Hence, $|N_{G'}(v)\cap \bar
S|\leq \delta_G(v)-1$, while $|N_{G'}(v)\cap S|\geq
\delta_G(v)+(r-1)+1=\delta_G(v)+r$. Therefore, $S$ is a valid
$r$-GOA.

Conversely, let $S$ be a $r$-GOA of $G'$. Since $S$ is a dominating
set, for each $K_2$-copy attached to $G$, either the corresponding
$A$- or the corresponding $B$-vertex is in $S$. Consider some $v\in
V\setminus S$. $v$ must be dominated. If no neighbor of $v$ in $V$
is in $S$, then $|N_{G'}(v)\cap S|\leq \delta_G(v)+r-1$, while
$|N_{G'}\cap \bar S|\geq \delta_G(v)$, which leads to a
contradiction. Hence, $S\cap V$ is a dominating set in $G$.

Combining the arguments, we obtain: $G=(V,E)$ has a dominating set
of size at most $k$ iff $G'=(V',E')$ has a $r$-GOA of size
$k+\sum_v(\delta_G(v)+r-1)= k+(r-1)|V|+2|E|$.

Now, we consider the case $r\geq 2$. Let $(G,k)$ be an instance of
Dominating Set with minimum degree $1$, with $G=(V,E)$. To any $v\in
V$, attach $\delta_G(v)+r-1\geq 1$ so-called $A$-vertices. All
$A$-vertices together form an independent set. Let
$A(v)=\{(v,1),\dots,(v,\delta_G(v)+r-1)\}$ denote the set of
$A$-vertices attached to $v\in V$. We denote the $B$-vertices
attached to the $A$-vertices in $A(v)$ by $B(v)$ and can describe
them as $B(v)=\left(
\begin{array}{c}A(v)\\r\end{array}\right)$, i.e., the
$r$-element subsets of $A(v)$. Each $X\in B(v)$ has as neighbors
exactly the $A$-vertices listed in $X$. This describes the graph
$G'=(V',E')$ as obtained from $G$.

If $D\subseteq V$ is a dominating set in $G$, then $S=D\cup A$ is a
$r$-GOA in $G'$. Clearly, $S$ is a dominating set in $G'$. Now,
consider a $B$-vertex $v$. Obviously, $N(v)\subseteq A(v)$, and
therefore $|N_{G'}(v)\cap S|=r\geq |N_{G'}(v)\cap\bar S|+r$. Any
vertex  $v\in V\setminus D$ has a neighbor $d\in D$. Hence,
$|N_{G'}(v)\cap \bar S|\leq \delta_G(v)-1$, while $|N_{G'}(v)\cap
S|\geq \delta_G(v)+(r-1)+1=\delta_G(v)+r$. Therefore, $S$ is a valid
$r$-GOA.

Conversely, let $S$ be a $r$-GOA of $G'$ of size $k+|A|$. Notice
that this bound is met if $S\cap V$ is a dominating set in $G$ and
all $A$-vertices go into $S$. Consider an $A(v)$-vertex $x$ and
assume $x\notin S$. Then, either there is a $y\in S\cap N(x)\cap
B(v)$, or $v\in S$, since otherwise $x$ would not be dominated.
Altogether, $x$ has $\left(
\begin{array}{c}\delta_G(v)+r-1\\r\end{array}\right)+1$
many neighbors. Since $S$ is an $r$-GOA, more than
$|A(v)|=\delta_G(v)+r-1$ vertices from the gadget attached to $v$
would be in $S$, this way violating the bound on the size of $S$.
 Consider some $v\in V\setminus S$. $v$ must be dominated.
If no neighbor of $v$ in $V$ is in $S$, then $|N_{G'}(v)\cap S|\leq
\delta_G(v)+r-1$, while $|N_{G'}\cap \bar S|\geq \delta_G(v)$, which
leads to a contradiction. Hence, $S\cap V$ is a dominating set in
$G$.

Combining the arguments, we obtain: $G=(V,E)$ has a dominating set
of size at most $k$ iff $G'=(V',E')$ has a $r$-GOA of size
$k+\sum_v(\delta_G(v)+r-1)= k+(r-1)|V|+2|E|$.
\end{proof}

\section{Bounding the offensive $r$-alliance number}

\begin{theorem}\label{cotainfsup1000}
For any graph $G$  of order  $n$ and minimum degree $\delta$, and
for every  $r\in \{2-\delta,\dots, \delta\}$,
$$\left\lceil\frac{\delta+r}{2}\right\rceil \le
a_{r}^{o}(G)\le \gamma_r^o(G) \le
n-\left\lceil\frac{\delta-r+2}{2}\right\rceil.$$
\end{theorem}

\begin{proof}
Let $v$ be a vertex of minimum degree in $G$ and let $Y\subset
N_{V}(v)$ such that $|Y|=\left\lceil\frac{\delta+r}{2}\right\rceil$.
Let $S=\{v\}\cup N_V(v)-Y$. Hence, $\bar{S}$ is a dominating set and
$$\delta_{\bar{S}}(v)=\left\lceil\frac{\delta+r}{2}\right\rceil \ge
\left\lfloor\frac{\delta+r}{2}\right\rfloor
=\delta-\left\lceil\frac{\delta+r}{2}\right\rceil+r =
\delta_{S}(v)+r.$$ Thus,
$$\delta_{\bar{S}}(u)\ge \delta_{\bar{S}}(v) \ge \delta_S(v)+r \ge \delta_S(u)+r, \quad \forall u\in S.$$
Therefore, $\bar{S}$ is a global offensive $r$-alliance in $G$ and,
as a consequence, the upper bound follows.

On the other hand, let  $X\subset V$ be an offensive $r$-alliance in
$G$. For every $v\in \partial (X)$ we have
  $$\delta(v) =\delta_X(v)+\delta_{\bar{X}}(v)$$
         $$ \delta(v)  \le \delta_X(v)+ \frac{\delta(v)-r}{2}$$
           $$ \frac{\delta(v)+r}{2} \le \delta_X(v)\le |X|$$
         $$ \frac{\delta+r}{2}\le |X|.$$
Therefore, the lower bound follows.
\end{proof}

The bounds are attained for every $r$ in the case of the complete
graph $G=K_n$.


A set $S\subset V$ is a $k$-dominating set if for every $v\in
\bar{S}$, $\delta_{S}(v)\geq k$. The \emph{ $k$-domination number}
of $G$, $\gamma_{_k}(G)$, is the minimum cardinality of a
$k$-dominating set in  $G$. The following result generalizes, to $r$
alliances, some previous results obtained for $r=1$ and $r=2$
\cite{albesiga2,albesiga7}.

\begin{theorem}\label{thOf}
For any simple graph $G$ of order $n$, minimum degree $\delta$, and
Laplacian spectral radius\footnote{i.e., the largest Laplacian eigenvalue
of $G$. The reader is referred to \cite{cv,Mohar91} for a detailed study and survey on
the Laplacian matrix of  a graph and its eigenvalues.}  $\mu_{*}$,
$$ \left\lceil\frac{n}{\mu_{*}}\left\lceil\frac{\delta
+r}{2}\right\rceil\right\rceil \le \gamma_{r}^{o}(G) \le
\left\lfloor\displaystyle\frac{\gamma_{_r}(G) +
n}{2}\right\rfloor.$$
\end{theorem}

\begin{proof}
Let $H\subset V$ be an $r$-dominating set of $G$ of minimum
cardinality. If $|\bar{H}|=1$, then $\gamma_r(G)=n-1$ and
$\gamma_{r}^{o}(G)\le n-1$. If $|\bar{H}|\neq 1$, let $\bar{H}=X\cup
Y$ be a partition of  $\bar{H}$ such that the edge-cut between $X$
and $Y$ has the maximum cardinality. Suppose $|X|\le |Y|$. For every
$v\in Y$, $\delta_{H}(v)\ge r$ and $\delta_{X}(v)\ge \delta_{Y}(v)$.
Therefore, the set  $W=H\cup X$ is a global offensive $r$-alliance
in $G$, i.e., for every $v\in Y$, $\delta_{W}(v)\ge
\delta_{Y}(v)+r$. Then we have,
\begin{equation}\label{eq3} 2|X|+\gamma_r(G)\le n\end{equation}
and
\begin{equation}\label{eq4} \gamma_{r}^{o}(G)\le|X|+\gamma_{r}(G) .\end{equation}
 Thus, by (\ref{eq3}) and (\ref{eq4}), we obtain the upper bound.

It was shown in  \cite{fiedler} that the Laplacian spectral radius
of $G$, $\mu_{*}$, satisfies
\begin{equation}\label{rfiedler2}
  \mu_{*}=2n \max \left\{ \frac{\displaystyle{\sum_{v_i\sim v_j}}(w_i-w_j)^2 }
  {\displaystyle{\sum_{v_i\in V}\sum_{v_j\in V}}(w_i-w_j)^2}: \mbox{\rm $w\neq \alpha{\bf j}$
   for  $\alpha\in \mathbb{R}$ } \right\},
\end{equation}
where $V=\{v_1, v_2, ..., v_n\}$, ${\bf j}=(1,1,...,1)$ and $w\in
\mathbb{R}^n$. Let $S\subset V$.  From  (\ref{rfiedler2}), taking
$w\in \mathbb{R}^n$ defined as
$$w_i= \left\lbrace \begin{array}{ll} 1  & {\rm if }\quad  v_i\in S;
                            \\ 0 &  {\rm otherwise} \end{array}
                                                \right .$$
we obtain
\begin{equation}\label{fiedlerAlliance100}
\mu_{*}\ge \frac{n\displaystyle\sum_{v\in \bar{S}} \delta
_S(v)}{|S|(n-|S|)}.
\end{equation}
Moreover, if $S$ is a global offensive $r$-alliance in  $G$,
\begin{equation} \label{globalGrado100}
\delta _S(v)\ge \left\lceil\frac{\delta(v)+r}{2}\right\rceil, \quad
\forall v\in \bar{S}.
\end{equation}
Thus, (\ref{fiedlerAlliance100}) and (\ref{globalGrado100}) lead to
\begin{equation}\label{final100}
\mu_{*}\ge \frac{n}{|S|}\left\lceil\frac{ \delta +r}{2}\right\rceil.
\end{equation}
Therefore, solving (\ref{final100}) for $|S|$ we obtain the lower
bound.\end{proof}

The above-mentioned bounds are attained, for instance, in the case of the
complete graph of order $n$.

\begin{corollary} For any simple graph $G$ of order $n$, minimum degree $\delta$,
and for every $r \in \{1,...,  \delta\}$, $$\gamma_{r}^{o}(G)\le
 \left\lfloor\displaystyle\frac{n(2r+1)}{2r+2}\right\rfloor.$$
\end{corollary}

\begin{proof}
The bound immediately  follows from the following bound on
$\gamma_r(G)$ \cite{Cockayne}:
\begin{equation}\label{eq5} \delta\ge
 r\Rightarrow\gamma_r(G)\le \frac{rn}{r+1}.
\end{equation}
\end{proof}

\begin{corollary} \label{coro21nue}
Let ${\cal L}(G)$ be the line graph of a $\delta$-regular graph $G$
of order $n$. Then
$$\gamma_{r}^o({\cal L}(G)) \ge \frac{n}{4}\left\lceil \frac{2(\delta-1)+r}{2}\right\rceil .$$
\end{corollary}

\begin{proof}
We denote by $A$ the adjacency matrix of ${\cal L}(G)$ and by
$2(\delta-1)=\lambda_0>\lambda_1>\cdots>\lambda_b=-2$ its distinct
eigenvalues. We denote by $L$ the Laplacian matrix of ${\cal
L}(\Gamma)$ and by $\mu_0=0<\mu_1<\cdots <\mu_b$ its distinct
Laplacian eigenvalues. Then, since $L= 2(\delta-1)I_n- A$, the
eigenvalues of both matrices, $A$ and $L$, are related by
\begin{equation} \label{eigen}
\mu_l = 2(\delta-1)-\lambda_l, \quad l=0, \dots, b.
\end{equation}
Thus, the Laplacian spectral radius of ${\cal L}(G)$ is
$\mu_b=2\delta$. Therefore, the result immediately follows.
\end{proof}

There are some immediate bounds on $\gamma_{_r}^o(G)$ derived from
the following remarks.

\begin{remark}
If $S$ is an independent set in $G$, then $\bar{S}$ is a global
offensive $r$-alliance in $G$ ($r\le  \delta$).
\end{remark}

\begin{remark}
All global offensive $r$-alliance in $G$ is a $\left\lceil
\frac{\delta+r}{2}\right\rceil$-dominating set in $G$ ($r\ge
2-\delta$).
\end{remark}

Therefore, the following bounds follow.
\begin{equation}\gamma_{_{\left\lceil \frac{\delta+r}{2} \right\rceil}}(G)\le
\gamma_{_r}^o(G)\le n-\alpha(G),
\end{equation}
 where $\alpha(G)$ denotes the independence number of $G$.

The reader is referred to our previous works
\cite{albesiga2,albesiga3,albesiga6,albesiga7} for a more detailed
study on offensive 1-alliances and offensive 2-alliances.

\end{document}